\documentclass[12pt]{amsproc}%{amsart}
\usepackage{amsmath,amssymb,amscd}

\begin{document}

\title[Free Probability for Pairs of Faces III]{Free Probability for Pairs of Faces III: \\
$2$-Variables Bi-free Partial $S$- and $T$-Transforms}
\author{Dan-Virgil Voiculescu}
\address{D.V. Voiculescu \\ Department of Mathematics \\ University of California at Berkeley \\ Berkeley, CA\ \ 94720-3840}
\thanks{Research supported in part by NSF Grant DMS-1301727.}
\date{}

\begin{abstract}
We introduce two $2$-variables transforms: the partial bi-free $S$-transform and the partial bi-free $T$-transform. These transforms are the analogues for the bi-multiplicative and, respectively, for the additive-multiplicative bi-free convolution of the $2$-variables partial bi-free $R$-transform in our previous paper in this series.
\end{abstract}

\maketitle

\setcounter{section}{-1}
\section{Introduction}
\label{sec0}

In the first paper in this series \cite{9} we proposed an extension of free probability to systems with left and right variables, based on the notion of bi-freeness. Due to the recent work (\cite{10}, \cite{7}, \cite{1}, \cite{2}, \cite{8}, \cite{3}, \cite{5}) we already have a much better understanding of bi-free probability.

Here, we complement our paper \cite{10} on the $2$-variables bi-free $R$-transform with two further $2$-variables partial transforms adapted to the bi-multiplicative operation $\boxtimes\boxtimes$ and to the additive-multiplicative operation $\boxplus\boxtimes$. In view of free-probability it was obvious to use the letter $S$ for the first transform, while for the second which corresponds to passing from bi-free pairs $(a_1,b_1)$, $(a_2,b_2)$ to $(a_1+a_2,b_1b_2)$ we went for $T$ the next letter in the alphabet.

Like in the case of addition \cite{10} we found that also for multiplication instead of our original approach to the one-variable $S$-transform \cite{12}, one of the proofs of U.~Haagerup in \cite{6} was very suitable as a starting point for the bi-free generalization. The additive and multiplicative cases being closely related in \cite{6}, we will be able to make also substantial use here of our work in \cite{10} on the bi-free partial $R$-transform. Since a combinatorial proof for the partial $R$-transform was found in \cite{8} and there is also combinatorial work for handling multiplication of bi-free systems of variables using additive bi-free cumulants \cite{1} it is to be expected that a combinatorial approach to the $2$-variables partial $S$- and $T$-transforms we consider here is possible. On the other hand, the study of infinite divisibility for additive bi-free convolution of probability measures on ${\mathbb R}^2$ in \cite{5} may not have an immediate analogue for the operations considered here, since the result of the operations $\boxtimes\boxtimes$ and $\boxplus\boxtimes$ on probability measures seems to produce signed measures only in general.

The paper has five sections including the introduction and the preliminaries in section~1. The next two sections, sections~2 and 3, deal with the essential parts of the proofs for the partial $S$- and respectively $T$-transforms. In the last section, section~4, we collect the immediate consequences of the work in the previous two sections.

\section{Preliminaries}
\label{sec1}

\subsection{Bi-freeness.}
\label{sec1.1} For the definition of when two two-faced pairs $(a_1,b_1)$ and $(a_2,b_2)$ in a non-commutative probability space $({\mathcal A},\varphi)$ are bi-free we refer the reader to \cite{9}. We should note that in the present paper we actually will not have to take our arguments all the way to the technical definition of bi-freeness, since we will in essence take as starting point a result in \cite{10}, which is Lemma~2.3 of the present paper, giving a certain consequence of bi-freeness. We will also use some of the simplest properties of bi-freeness and the definition of multiplicative bi-free convolution \cite{9}.

\subsection{Notations.}
\label{sec1.2} The notation we will use like in \cite{10} is a combination of the usual notation from free probability for the $S$-transform (\cite{12},\cite{13}) with notations used in the paper of Haagerup (\cite{6}) and some natural two-variable extensions of these.

If $({\mathcal A},\varphi)$ is a Banach-algebraic non-commutative probability space, that is, ${\mathcal A}$ is a unital Banach algebra over ${\mathbb C}$ and $\varphi: {\mathcal A} \to {\mathbb C}$ is a continuous functional so that $\varphi(1) = 1$. If $a \in {\mathcal A}$ and $t,z \in {\mathbb C}$ are so that $z1-a$ and $1-ta$ are invertible we shall consider $G_a(z) = \varphi((z1-a)^{-1})$ and $h_a(t) = \varphi((1-ta)^{-1})$. If $t \ne 0$, then $h_a(t) = t^{-1}G_a(t^{-1})$. In the definition of the $S$-transform we also use $\psi_a(t) = h_a(t) - 1$, which in part of \cite{6} is also denote $l_a(t)$. We shall also deal with certain two variables analogues $G_{a,b}(z,w) = \varphi((z1-a)^{-1}(w1-b)^{-1})$ and $H_{a,b}(t,s) = \varphi((1-ta)^{-1}(1-sb)^{-1})$ where $a,b \in {\mathcal A}$, $z,w,s,t \in {\mathbb C}$ are such that the inverses of $z1-a,w1-b,1-ta,1-sb$ exist.

If $\varphi(a) \ne 0$, then $\psi_a(0) = 0$, $\psi'_a(0) = \varphi(a) \ne 0$ so that $\psi_a$ has an inverse function ${\mathcal X}_a(z)$ in a neighborhood of $0$. The $S$-transform is then $S_a(z) = \frac {1+z}{z} {\mathcal X}_a(z)$. The key property of $S_a(z)$ is that if $a_1,a_2 \in ({\mathcal A},\varphi)$ with $\varphi(a_1) \ne 0$, $\varphi(a_2) \ne 0$ independent, then $S_{a,a_2}(z) = S_{a_1}(z)S_{a_2}(z)$.

We also recall the $R$-transform (\cite{11}, \cite{13}). Let $K_a(z)$ be the inverse of $G_a(z)$ in a neighborhood of $0$ and then $R_a(z) = K_a(z)-z^{-1}$.

We will be dealing with these quantities as germs of holomorphic functions and in this vein say that $h_a(t)$ is ``holomorphic near $0$'' or $H_{a,b}(t,0)$ ``holomorphic near $(0,0) \in {\mathbb C}^2$'', etc.

We should remark that frequently the indices for previous quantities are $\mu_a,\mu_b,\mu_{a,b}$ that is the distributions of the variables, instead of the variables themselves. Clearly, the quantities
\[
G_a,h_a,\psi_a,G_{a,b},H_{a,b},{\mathcal X}_a,S_a
\]
depend only on the distributions of the variables, the present notation using the variables is just a convenient abbreviation.

\subsection{Formal Power Series.}
\label{sec1.3} The quantities
\[
G_a,h_a,\psi_a,G_{a,b},H_{a,b},{\mathcal X}_a,S_a
\]
can also be defined in the purely algebraic context when ${\mathcal A}$ is an algebra over ${\mathbb C}$ and $\varphi: {\mathcal A} \to {\mathbb C}$ is a linear map.  Then
\[
\begin{split}
G_a(z) &= \sum_{n \ge 0} z^{-n-1}\varphi(a^n), \\
h_a(t) &= \sum_{k \ge 0} t^k\varphi(a^k), \\
G_{a,b}(z,w) &= \sum_{m \ge 0,n \ge 0} z^{-m-1}w^{-n-1}\varphi(a^mb^n), \\
H_{a,b}(t,s) &= \sum_{p \ge 0,q \ge 0} t^ps^q\varphi(a^pb^q),
\end{split}
\]
and if $\varphi(a) \ne 0$ also ${\mathcal X}_a(z)$ and $S_a(z) = \sum_{k ge 0} s_kz^k$ are defined as formal power series. See $4.5$ below concerning why proofs carried out in the Banach algebra context in the end give results which also hold for formal power series.

\clearpage
\section{The partial bi-free $S$-transform $S_{a,b}(z,w)$ for a two-faced pair of variables}
\label{sec2}

\subsection{}
\label{sec2.1} Throughout section~\ref{sec2} $({\mathcal A},\varphi)$ will be a unital Banach algebra ${\mathcal A}$ endowed with a unit-preserving linear expectation functional $\varphi$.

\bigskip
\noindent
{\bf Definition 2.1.} If $(a,b)$ is a two-faced pair of non-commutative random variables in $({\mathcal A},\varphi)$ so that $\varphi(a) \ne 0$, $\varphi(b) \ne 0$, we define the $2$-variables partial bi-free $S$-transform as the holomorphic function of $(z,w) \in ({\mathbb C}\backslash\{0\})^2$ when $z,w$ are near $0$,
\[
S_{a,b}(z,w) = \frac {z+1}{z} \frac {w+1}{w} \left( 1 - \frac {1+z+w}{H_{a,b}({\mathcal X}_a(z),{\mathcal X}_b(w))}\right).
\]
Since this involves only the joint distribution $\mu_{a,b}$ of $(a,b)$, we shall also write $S_{\mu_{a,b}}(z,w)$ for $S_{a,b}(z,w)$.

\bigskip
Our aim in section~\ref{sec2} will be to prove the following theorem.

\bigskip
\noindent
{\bf Theorem 2.1.} {\em Let $(a_1,b_1)$ and $(a_2,b_2)$ be a bi-free pair of two-faced pairs of non-commutative random variables in $({\mathcal A},\varphi)$ so that $\varphi(a_1) \ne 0$, $\varphi(a_2) \ne 0$, $\varphi(b_1) \ne 0$, $\varphi(b_2) \ne 0$. Then we have
\[
S_{a_1a_2,b_1b_2}(z,w) = S_{a_1,b_1}(z,w)S_{a_2,b_2}(z,w)
\]
for $(z,w) \in ({\mathbb C}\backslash\{0\})^2$ near $(0,0)$.
}

\subsection{}
\label{sec2.2} In preparation for the proof of the theorem we will have a few lemmas. The following result is from the paper by Haagerup (\cite{6}), a proof is given for the reader's convenience.

\bigskip
\noindent
{\bf Lemma 2.1.} {\em Let $a_1,a_2 \in ({\mathcal A},\varphi)$. If $t_1,t_2 \in {\mathbb C}\backslash\{0\}$ are near zero and satisfy
\[
h_{a_1}(t_1) = h_{a_2}(t_2) = h \notin \{0,1\}
\]
then with $\rho = (h(h-1))^{-1}$ and
\[
t = \frac {h}{h-1} t_1t_2
\]
we have
\[
(1-t_1a_1)(1-\rho a_1(t_1)a_2(t_2))(1-t_2a_2) = \frac {1}{h} (1-ta_1a_2).
\]
}

\bigskip
\noindent
{\bf {\em Proof.}} We have expanding
\[
(1-t_1a_1)(1-\rho a_1(t_1)a_2(t_2))(1-t_2a_2) = c_0+c_1a_1+c_2a_2+c_3a_1a_2
\]
where
\[
\begin{split}
c_3 &= t_1t_2(1-\rho h_{a_1}(t_1)h_{a_2}(t_2)) \\
\  &= t_1t_2(1 - (h(h-1))^{-1}h^2) = -\frac {t_1t_2}{h-1}\,, \\
c_0 &= 1 - \rho(1-h_{a_1}(t_1))(1 - h_{a_2}(t_2)) \\
\  &= 1 - \frac {(1-h)^2}{h(h-1)} = \frac {1}{h}\,, \\
c_1 &= -t_1(1+\rho h_{a_1}(t_1) - \rho h_{a_1}(t_1)h_{a_2})(t_2)) \\
\  &= -t_1\left( 1 + \frac {h-h^2}{h(h-1)}\right) = 0, \\
c_2 &= -t_2(1-1) = 0,
\end{split}
\]
similarly. Thus $c_1 = c_2 = 0$ and
\[
c_3/c_0 = -t_1t_2 \frac {h}{h-1} = -t.
\]
This gives
\[
\begin{split}
c_0+c_1a_1+c_2a_2+c_3a_1a_2 &= c_0\left( 1 + \frac {c_3}{c_0} t_1t_2\right) \\
&= \frac {1}{h} (1 - ta_1a_2).
\end{split}
\]
\qed

\subsection{}
\label{sec2.3} In case $a_1,a_2$ are freely independent we can produce $t_1,t_2$ satisfying the assumptions of the preceding lemma using the quantities in the definition of the $S$-transform (\cite{12}).

\bigskip
\noindent
{\bf Lemma 2.2.} {\em Assume $a_1,a_2 \in ({\mathcal A},\varphi)$ are freely independent and $\varphi(a_k) \ne 0$, $k = 1,2$. Then for $z$ sufficiently close to $0$ taking $t_k = {\mathcal X}_{a_k}(z)$, $k = 1,2$, we have that
\[
h_{a_1}(t_1) = h_{a_2}(t_2) = z + 1
\]
and if, moreover, $z \ne 0$ then we also have
\[
t = \frac {z+1}{z} t_1t_2 = {\mathcal X}_{a_1a_2}(z) \text{ and } h_{a_1a_2}(t) = z+1.
\]
Further, if $z \ne 0$ and $\rho = (z(1+z))^{-1}$, then
\[
(1-t_1a_1)(1-\rho a_1(t_1)a_2(t_2))(1-t_2a_2) = \frac {1}{z+1}(1-ta_1a_2).
\]
}

\bigskip
\noindent
{\bf {\em Proof.}} Since $\psi_{a_k}(z) = h_{a_k}(z) - z$ and ${\mathcal X}_{a_k}$ is the inverse of $\psi_{a_k}$ we get that $t_k = {\mathcal X}_{a_k}(z)$ will satisfy $h_{a_k}(t_k) = z+1$, $k = 1,2$, and we can apply Lemma~2.1 to these $t_1,t_2$ to get the remaining assertions. Note that if $z \ne 0$ is sufficiently close to $0$, then $z + 1 \notin \{0,1\}$. Then also
\[
\begin{split}
t &= \frac {z+1}{z} t_1t_2 = \frac {z}{z+1} \left( \frac {z+1}{z} t_1\right) \left( \frac {z+1}{z} t_2\right) \\
&= \frac {z}{z+1} S_{a_1}(z)S_{a_2}(z) = \frac {z}{z+1} S_{a_1a_2}(z) = {\mathcal X}_{a_1a_2}(z).
\end{split}
\]
\qed

\subsection{}
\label{sec2.4} The next result we record as a lemma is a variant of a result from Part~II (\cite{10}, Lemma~2.9). The proof being the same, with the only difference that $\|\rho a_1(t_1)a_2(t_2)\| < 1$, $\|\sigma b_1(s_1)b_2(s_2)\| < 1$ are now among the assumptions, instead of assuming $|\rho| < C$, $|\sigma| < C$ and deriving these inequalities from $\|a_k(t_k)\| \to 0$, $\|b_k(s_k)\| \to 0$ as $t_k \to 0$, $s_k \to 0$.

\bigskip
\noindent
{\bf Lemma 2.3.} {\em Assume the two-faced pairs $(a_1,b_1)$, $(a_2,b_2)$ in $({\mathcal A}
,\varphi)$ are bi-free and that $\|\rho a_1(t_1)a_2(t_2)\| < 1$, $\|\sigma b_1(s_1)b_2(s_2)\| < 1$, $\|t_ka_k\| < 1$, $\|s_kb_k\| < 1$, where $\rho,\sigma,t_1,t_2,s_1,s_2 \in {\mathbb C}$. Then we have
\[
\begin{split}
&\varphi(((1-t_1a_1)(1-\rho a_1(t_1)a_2(t_2))(1-t_2a_2))^{-1} \\
&\hskip 2 in ((1-s_1b_1)(1-\sigma b_1(s_1)b_2(s_2))(1-s_2b_2))^{-1}) \\
&= \frac {(\varphi((a_1(t_1)b_1(s_1)) + h_{a_1}(t_1)h_{b_1}(s_1))(\varphi((a_2(t_2)b_2(s_2)) + h_{a_2}(t_2)h_{b_2}(s_2))))}{1-\rho\sigma\varphi(a_1(t_1)b_1(s_1))\varphi(a_2(t_2)b_2(s_2))}\,.
\end{split}
\]
}

\subsection{}
\label{sec2.5} We can now start combining the preceding results in the next lemma.

\bigskip
\noindent
{\bf Lemma 2.4.} {\em Assume the two-faced pairs $(a_1,b_1)$ and $(a_2,b_2)$ in $({\mathcal A},\varphi)$ are bi-free and assume that $\varphi(a_k) \ne 0$, $\varphi(b_k) \ne 0$, $k = 1,2$. If $z,w \in {\mathbb C}\backslash\{0\}$ are sufficiently close to $0$ then
\[
t = {\mathcal X}_{a_1a_2}(z),\ t_k = {\mathcal X}_{a_k}(z),\ s = {\mathcal X}_{b_1b_2}(w),
\]
$s_k = {\mathcal X}_{b_k}(z)$, $k = 1,2$, will be $\ne 0$ and will also satisfy $\|(z(z+1))^{-1}a_1(t_1)a_2(t_2)\| < 1$ and $\|(w(w+1))^{-1}b_1(s_1)b_2(s_2)\| < 1$ and
\[
\begin{split}
h_{a_1a_2}(t) &= h_{a_1}(t_1) = h_{a_2}(t_2) = z + 1 \notin \{0,1\} \\
h_{b_1b_2}(s) &= h_{b_1}(s_1) = h_{b_2}(s_2) = w + 1 \notin \{0,1\}.
\end{split}
\]
Moreover, we have:
\[
\begin{split}
&(z+1)(w+1)H_{a_1a_2,b_1b_2}(t,s) \\
&= \frac {H_{a_1,b_1}(t_1,s_1)H_{a_2,b_2}(t_2,s_2)}{1 - \frac {1}{zw(z+1)(w+1)} (H_{z_1,b_1}(t_1,s_1)-(z+1)(w+1))(H_{a_2,b_2}(t_2,s_2) - (z+1)(w+1))}\,.
\end{split}
\]
}

\bigskip
\noindent
{\bf {\em Proof.}} As we saw in Lemma~$2.2$, since $\psi_a(z) = h_a(z) - 1$, if $\varphi(a) \ne 0$, so that $\psi'_a(0) = \varphi(a) \ne 0$ and ${\mathcal X}_a(z)$ exists, we have $h_a({\mathcal X}_a(z)) = z+1$ for $z$ in a neighborhood of $0$. This takes care of the fact that the assumptions for the application of Lemma~$2.1$ to $a_1,a_2,t_1,t_2$ and, respectively, $b_1,b_2,s_1,s_2$ are satisfied. The fact that $\|(z(z+1))^{-1}a_1(t_1)a_2(t_2)\| \to 0$ as $z \to 0$, $z \ne 0$ follows from
$\|a_k(t_k)\| = O(|t_k|) = O(|{\mathcal X}_{a_k}(z)|) = O(|z|)$ so that $\|(z(z+1))^{-1}a_1(t_1)a_2(t_2)\| = O(|z|)$ as $z \to 0$. Similarly, if $w \ne 0$, $\|(w(w+1))^{-1}b_1(s_1)b_2(s_2)\| \to 0$ as $w \to 0$. This verifies that the assumptions for the application of Lemma~2.3 are satisfied.

In view of Lemma~2.1, the left-hand side of the equality in Lemma~2.3 becomes for the present choice of $t_1,t_2$
\[
(z+1)(w+1)\varphi((1-ta_1a_2)^{-1}(1-sb_1b_2)^{-1}) = (z+1)(w+1)H_{a_1a_2,b_1b_2}(t,s).
\]

On the other hand, the numerator of the right-hand side of the equality in Lemma~2.3 is
\[
H_{a_1,b_1}(t_1,s_1)H_{a_2,b_2}(t_2,s_2)
\]
because 
\[
\begin{split}
H_{a_1,b_1}(t_1,s_1) &= \varphi((a_1(t_1)+h_{a_1}(t_1)1)(b_1(s_1)+h_{b_1}(s_1)1)) \\
&= \varphi(a_1(t_1)b_1(s_1)) + h_{a_1}(t_1)h_{b_2}(s_2)
\end{split}
\]
and a similar equality holds for $H_{a_2,b_2}(t_2,s_2)$. Further on, the denominator of the right-hand side in Lemma~2.3 is
\[
\begin{split}
&1 - \rho\sigma\varphi(a_1(t_1)b_1(s_1))\varphi(a_2(t_2)b_2(s_2)) = \\
&1 - \frac {1}{z(z+1)w(w+1)} (H_{a_1,b_1}(t_1,s_1) \\
&\ \  - (z+1)(w+1))(H_{a_2,b_2}(t_2,s_2) - (z+1)(w+1))
\end{split}
\]
which concludes the proof.\qed

\subsection{Proof of Theorem 2.1.}
\label{sec2.6} The equality to be proved in Theorem~2.1 is just another form of the equality in Lemma~2.4 and can be obtained from the lemma by some simple algebraic manipulation.

It will be convenient to use some abbreviated notation:
\[
\begin{split}
&(z+1)(w+1) = \beta,\ zw = \alpha \\
&H_{a_1,b_1}(t_1,s_1) = H_1,\ H_{a_2,b_2}(t_2,s_2) = H_2 \\
&H_{a_1a_2,b_1b_2}(t,s) = H_{12}.
\end{split}
\]
Since we have $t_1 = {\mathcal X}_{a_1}(z)$, $t_2 = {\mathcal X}_{a_2}(z)$, $t_3 = {\mathcal X}_{a_3a_2}(z)$ and similar formulas for $s_1,s_2,s$, we see that the equality to be proved in Theorem~2.1 is:
\[
\frac {\beta}{\alpha} \left( 1 + \frac {\alpha-\beta}{H_{12}}\right) = \frac {\beta}{\alpha} \left( 1 + \frac {\alpha-\beta}{H_1}\right) \frac {\beta}{\alpha} \left( 1 + \frac {\alpha-\beta}{H_2}\right)\,.
\]
On the other hand the equality in the lemma after taking inverses gives
\[
\frac {1}{\beta H_{12}} = \frac {1}{H_1H_2} - \frac {1}{\alpha\beta} \frac {(H_1-\beta)(H_2-\beta)}{H_1H_2}\,.
\]
This gives
\[
\frac {1}{\beta H_{12}} = \frac {1}{H_1H_2} \left( 1 - \frac {\beta}{\alpha}\right) + \frac {1}{\alpha} \left( \frac {1}{H_1} + \frac {1}{H_2}\right) - \frac {1}{\alpha\beta}\,.
\]
The right-hand side becomes successively
\[
\begin{split}
&\frac {1}{\alpha(\alpha-\beta)} \left( \frac {(\alpha-\beta)^2}{H_1H_2} + \frac {\alpha-\beta}{H_1} + \frac {\alpha-\beta}{H_2}\right) - \frac {1}{\alpha\beta} \\
&= \frac {1}{\alpha(\alpha-\beta)} \left( \frac {\alpha-\beta}{H_1} + 1\right) \left( \frac {\alpha-\beta}{H_2} + 1\right) - \frac {1}{\alpha} \left( \frac {1}{\alpha-\beta} + \frac {1}{\beta}\right) \\
&= \frac {1}{\alpha(\alpha-\beta)} \left( \frac {\alpha-\beta}{H_1} + 1\right) \left( \frac {\alpha-\beta}{H_2} + 1\right) - \frac {1}{\beta(\alpha-\beta)}\,.
\end{split}
\]
Hence the equality in the lemma gives
\[
\frac {1}{\beta(\alpha-\beta)} \left( \frac {\alpha-\beta}{H_{12}} + 1\right) = \frac {1}{\alpha(\alpha-\beta)} \left( \frac {\alpha-\beta}{H_1} + 1\right) \left( \frac {\alpha-\beta}{H_2} + 1\right)\,.
\]
Multiplying with $\frac {\beta^2(\alpha-\beta)}{\alpha}$ we get
\[
\frac {\beta}{\alpha} \left( \frac {\alpha-\beta}{H_{12}} + 1\right) = \frac {\beta^2}{\alpha^2} \left( \frac {\alpha-\beta}{H_1} + 1\right) \left( \frac {\alpha-\beta}{H_2} + 1\right)\,,
\]
which concludes the proof.\qed

\section{The partial bi-free $T$-transform $T_{a,b}(z,w)$ for a two-faced pair of variables}
\label{sec3}

\subsection{}
\label{sec3.1} Throughout section~\ref{sec3} $({\mathcal A},\varphi)$ will be a unital Banach algebra ${\mathcal A}$ endowed with a unit-preserving linear expectation functional $\varphi$. The $T$-transform will be a transform adapted to the bi-free additive-multiplicative convolution, which could be denoted $\boxplus\boxtimes$. It corresponds to the operation on bi-free two-faced pairs of variables $(a_1,b_1)$, $(a_2,b_2)$ which yields $(a_1+a_2,b_1b_2)$.

\bigskip
\noindent
{\bf Definition 3.1.} If $(a,b)$ is a two-faced pair in $({\mathcal A},\varphi)$, so that $\varphi(b) \ne 0$, we define the $2$-variables $T$-transform as the holomorphic function of $(z,w) \in ({\mathbb C}^2\backslash\{0\})^2$ near $(0,0) \in {\mathbb C}^2$
\[
T_{a,b}(z,w) = \frac {w+1}{w} \left( 1 - \frac {z}{F_{a,b}(K_a(z),{\mathcal X}_b(w))} \right)
\]
where $F_{a,b}(t,s) = \varphi((t-a)^{-1}(1-sb)^{-1})$. Since this involves only the joint distribution of $(a,b)$ we shall also write $T_{\mu_{a,b}}$ instead of $T_{a,b}$.

\bigskip
Our aim in section~\ref{sec3} will be to prove the following theorem.

\bigskip
\noindent
{\bf Theorem 3.1.} {\em Let $(a_1,b_1)$ and $(a_2,b_2)$ be a bi-free pair of two-faced pairs of non-commutative random variables in $({\mathcal A},\varphi)$, so that $\varphi(b_1) \ne 0$, $\varphi(b_2) \ne 0$. Then we have
\[
T_{a_1,b_1}(z,w)T_{a_2,b_2}(z,w) = T_{a_1+a_2,b_1b_2}(z,w)
\]
if $(z,w) \in ({\mathbb C}\backslash\{0\})^2$ are near $(0,0)$.
}

\subsection{}
\label{sec3.2} The proof will be along similar lines to the proofs of the partial bi-free $R$ and $S$ transforms. For the right half of $T$, which is multiplicative, we will be able to use lemmas from section~\ref{sec2}. To deal with the left half, the additive half, of $T_{a,b}$ some of the material used in \cite{10} can also be used here. Since this is another paper we will recall the facts we need, sometimes leaving out the proofs.

\subsection{}
\label{sec3.3} The lemma which follows gives a result from the paper by Haagerup \cite{6}, which we also used in \cite{10} and accompanied there by a proof, which will be omitted this time.

\bigskip
\noindent
{\bf Lemma 3.1.} {\em Let $a_1,a_2 \in ({\mathcal A},\varphi)$ and assume $t_1,t_2 \in {\mathbb C}$ are sufficiently close to zero and satisfy $t_1h_{a_1}(t_1) = t_2h_{a_2}(t_2)$ and $\rho = (h_{a_1}(t_1)h_{a_2}(t_2))^{-1}$ and $t = \frac {t_1h_{a_1}(t_1)}{h_{a_1}(t_1)+h_{a_2}(t_2)-1}$ are defined. Then we have
\[
(1-t_1a_1)(1-\rho a_1(t_1)a_2(t_2))(1-t_2a_2) = \frac {h_{a_1}(t_1)+h_{a_2}(t_2) - 1}{h_{a_1}(t_1)h_{a_2}(t_2)}\ (1 - t(a_1+a_2)).
\]
}

\subsection{}
\label{sec3.4} In case $a_1,a_2$ are freely independent we can produce $t_1,t_2$ satisfying the assumptions of the preceding lemma using the quantities involved in the definition of the $R$-transform (\cite{11}).

\bigskip
\noindent
{\bf Lemma 3.2.} {\em Assume $a_1,a_2$ are freely independent in $({\mathcal A},\varphi)$ then if $z \in {\mathbb C}\backslash\{0\}$ is sufficiently close to zero taking $t_k = (K_{a_k}(z))^{-1}$ we will have $t_1h_{a_1}(t_1) = t_2h_{a_2}(t_2) = z$ and
$\rho = (h_{a_1}(t_1)h_{a_2}(t_2))^{-1}$, $t = \frac {t_1h_{a_1}(t_1)}{h_{a_1}(t_1)+h_{a_2}(t_2)-1}$ will be defined and $|\rho| < 2$. Moreover, we will also have: $h_{a_1+a_2}(t) = h_{a_1}(t_1) + h_{a_2}(t_2) - 1$ and $th_{a_1+a_2}(t) = t_1h_{a_1}(t_1) = t_2h_{a_2}(t_2)$.
}

\bigskip
\noindent
{\bf {\em Proof.}} We have that if $t_k = (K_{a_k}(z))^{-1}$ then $t_kh_{a_k}(t_k) = G_{a_k}(t_k^{-1}) = G_{a_k}(K_{a_k}(z)) = z$. This also gives in view of the freeness of $a_1,a_2$ that if $z \ne 0$ is close to zero:
\[
\begin{split}
(K_{a_1+a_2}(z))^{-1} &= (K_{a_1}(z) + K_{a_2}(z) - z^{-1})^{-1} \\
&= (t_1^{-1} + t_2^{-1} - (h_{a_k}(t_k)t_k)^{-1})^{-1} \\
&= t_kh_{a_k}(t_k)(h_{a_1}(t_1) + h_{a_2}(t_2) - 1)^{-1} = t
\end{split}
\]
and hence $th_{a_1+a_2}(t) = G_{a_1+a_2}(K_{a_1+a_2}(z)) = z$. Moreover this also gives
\[
z(h_{a_1}(t_1) + h_{a_2}(t_2) - 1)^{-1} = t = z(h_{a_1+a_2}(t))^{-1}
\]
which implies
\[
h_{a_1+a_2}(t) = h_{a_1}(t_1) + h_{a_2}(t_2) - 1.
\]
Remark also that if $z \to 0$ then $K_{a_k}(z) \to \infty$ so that $t_k \to 0$ and $h_{a_k}(t_k) \to 1$. Thus if $z$ is close to $0$ we will have $\rho$ close to $1$, in particular $|\rho| < 2$.\qed

\subsection{}
\label{sec3.5} We will now apply Lemma~2.3 with $t_1,t_2,t,\rho$ chosen like in Lemma~3.2 and $s_1,s_2,s,\sigma$ as in Lemma~2.2.

\bigskip
\noindent
{\bf Lemma 3.3.} {\em Let $(a_1,b_1)$ and $(a_2,b_2)$ be two two-faced pairs in $({\mathcal A},\varphi)$, which are bi-free and assume $\varphi(b_k) \ne 0$, $k = 1,2$. Let further $z,w \in {\mathbb C}\backslash\{0\}$ be sufficiently close to $0$ and let $t = (K_{a_1+a_2}(z))^{-1}$, $t_k = (K_{a_k}(z))^{-1}$, $s = {\mathcal X}_{b_1b_2}(w)$, $s_k = {\mathcal X}_{b_k}(w)$, $k = 1,2$. Then we have:
\[
\begin{split}
&\frac {t}{t_1t_2} z(w+1)H_{a_1+a_2,b_1,b_2}(t,s) \\
&= \frac {H_{a_1,b_1}(t_1,s_1)H_{a_2,b_2}(t_2,s_2)}{1 - \frac {t_1t_2}{z^2} \cdot \frac {1}{w(w+1)}\left( H_{a_1,b_1}(t_1,s_1) - \frac {z(w+1)}{t_1}\right)\left( H_{a_2,b_2}(t_2,s_2) - \frac {z(w+1)}{2}\right)}\,.
\end{split}
\]
}

\bigskip
\noindent
{\bf {\em Proof.}} We shall apply Lemma~2.3 with $a_1,a_2,b_1,b_2,t_1,t_2,s_1,s_2$ as specified in the statement of the present lemma and $\rho = (h_{a_1}(t_1)h_{a_2}(t_2))^{-1}$, $\sigma = (h_{b_k}(s_k)(h_{b_k}(s_k)-1))^{-1}$. For $z$ close to $0$, $t_k$ are also close to $0$ and hence $\rho$ will be close to $1$ while $\|a_k(t_k)\|$ will also be close to $0$, so that the condition $\|\rho a_1(t_1)a_2(t_2)\| < 1$ will be satisfied. Similarly, if $w$ is close to $0$, as explained in the proof of Lemma~2.4, we have $\sigma = (w(w+1))^{-1}$ and we will have $\|\sigma b_1(s_1)b_2(s_2)\| < 1$ since it is $O(|w|)$ as $w \to 0$. Also clearly $\|t_ka_k\| < 1$, $\|s_kb_k\| < 1$ are satisfied for small $z$ and $w$.

In view of Lemma~3.2 and of Lemma~3.1 we have
\[
\begin{split}
&(1-t_1a_1)(1-\rho a_1(t_1)a_2(t_2))(1-t_2a_2) \\
&= \frac {h_{a_1}(t_1)+h_{a_2}(t_2)-1}{h_{a_1}(t_1)h_{a_1}(t_2)} (1 - t(a_1+a_2)) \\
&= \frac {h_{a_1+a_2}(t)}{h_{a_1}(t_1)h_{a_2}(t_2)} (1 - t(a_1+a_2)) \\
&= \frac {t_1t_2}{tz} (1 - t(a_1+a_2).
\end{split}
\]

Similarly using Lemma~2.2 and Lemma~2.1 we have
\[
(1 - s_1b_1)(1 - \sigma b_1(s_1)b_2(s_2))(1-s_2b_2) = \frac {1}{w+1} (1-sb_1b_2).
\]
It follows that the left-hand side of the formula in Lemma~2.3 equals
\[
\frac {tz(w+1)}{t_1t_2} H_{a_1+a_2,b_1b_2}(t,s).
\]
The numerator of the right-hand side in Lemma~2.3 always equals $H_{a_1,b_1}(t_1,s_1)H_{a_2,b_2}(t_2,s_2)$. To see that also the denominator in the right-hand side of the present lemma equals the denominator in the right-hand side of the formula in Lemma~2.3 it suffices to notice in addition to the fact that $\sigma = (w(w+1))^{-1}$ also that $\rho = (h_{a_1}(t_1)h_{a_2}(t_2))^{-1} = \frac {t_1t_2}{z^2}$ and that $h_{a_k}(t_k)h_{b_k}(s_k) = \frac {z}{t_k}(w+1)$.\qed

\subsection{Proof of Theorem~3.1}
\label{sec3.6} To prove the theorem we need to work more on the equality in Lemma~3.3. It will be convenient to abbreviate notation by writing $H_k$ for $H_{a_k,b_k}(t_k,s_k)$, $k = 1,2$ and $H_{12}$ for $H_{a_1+a_2,b_1b_2}(t,s)$. Taking inverses the equality in Lemma~3.3 becomes
\[
\frac {t_1t_2}{tz(w+1)} \cdot \frac {1}{H_{12}} = \frac {1}{H_1H_2} \left( 1 - \frac {w+1}{w}\right) - \frac {t_1t_2}{z^2w(w+1)} + \frac {1}{H_1} \cdot \frac {t_2}{zw} + \frac {1}{H_2} \cdot \frac {t_1}{zw}\,.
\]
We use now $tH_{a,b}(t,s) = F_{a,b}(t^{-1},s)$ and the abbreviations
\[
F_k = F_{a_k,b_k}(t^{-1}_k,s_k),\ F_{12} = F_{a_1+a_2,b_1b_2}(t^{-1},s).
\]
We have after multiplying by $(t_1t_2)^{-1}w$
\[
\begin{split}
\frac {w}{z(w+1)} \frac {1}{F_{12}} &= -\frac {1}{F_1F_2} - \frac {1}{z^2(w+1)} + \frac {1}{z} \cdot \frac {1}{F_1} + \frac {1}{z} \cdot \frac {1}{F_2} \\
&= -\left( \frac {1}{F_1} - \frac {1}{z}\right)\left( \frac {1}{F_2} - \frac {1}{z}\right) + \frac {w}{z^2(w+1)}\,.
\end{split}
\]
Multiplying with $-\left( \frac {z(w+1)}{w}\right)^2$ and moving terms we get
\[
\frac {w+1}{w} \left( 1 - \frac {z}{F_{12}}\right) = \frac {w+1}{w} \left( 1 - \frac {z}{F_1}\right) \cdot \frac {w+1}{w} \left( 1 - \frac {z}{F_2}\right)
\]
which is what we wanted to prove.\qed

\section{Concluding remarks}
\label{sec4}

\subsection{}
\label{sec4.1} After having established the transforms in the preceding two sections, it will be now easy to get some improvements of the results and derive some consequences in the present section.

\subsection{Analytic extension}
\label{sec4.2}

In Definitions~2.1 and 3.1 the $S$ and $T$ transforms were defined assuming that the variables satisfy $z \ne 0$ and $w \ne 0$. These conditions will be removed now.

\bigskip
\noindent
{\bf Proposition 4.1.} {\em 
\begin{itemize}
\item[a)] The transform $S_{a,k}(z,w)$ extends to a holomorphic function of $(z,w)$ in a neighborhood of $(0,0) \in {\mathbb C}^2$.
\item[b)] The transform $T_{a,b}(z,w)$ extends to a holomorphic function of $(z,w)$ ina neighborhood of $(0,0) \in {\mathbb C}^2$.
\end{itemize}
}

\bigskip
\noindent
{\bf {\em Proof.}}  a) We recall that if $\varphi(a) \ne 0$, $\varphi(b) \ne 0$
\[
S_{a,b}(z,w) = \frac {z+1}{z} \frac {w+1}{w} \left( 1 - \frac {1+z+w}{H_{a,b}({\mathcal X}_a(z),{\mathcal X}_b(w)}\right)
\]
was defined originally for $(z,w) \in ({\mathbb C}\backslash\{0\})^2$ close to $(0,0)$. We have $H_{a,b}(t,s) = \sum_{p \ge 0,q \ge 0} t^ps^q\varphi(a^pb^q)$ in a neighborhood of $(0,0)$ so that
\[
H_{a,b}(t,s) = h_a(t) + h_b(s) - 1 + st\eta(s,t)
\]
where $\eta(s,t)$ is a holomorphic function of $(s,t)$ in a neighborhood of $(0,0)$. Using $h_a({\mathcal X}_a(z)) = z + 1$, $h_b({\mathcal X}_b(w)) = w + 1$ we get
\[
S_{a,b}(z,w) = \frac {z+1}{z} \frac {w+1}{w} \frac {{\mathcal X}_a(z){\mathcal X}_b(w)\eta({\mathcal X}_a(z),{\mathcal X}_b(w))}{H_{a,b}({\mathcal X}_a(z),{\mathcal X}_b(w))}
\]
which gives a holomorphic function in a neighborhood of $(0,0)$ since $z^{-1}{\mathcal X}_a(z)$ and $w^{-1}{\mathcal X}_b(w)$ have removable singularities at $0$, their values at $0$ being $(\varphi(a))^{-1}$ and $(\varphi(b))^{-1}$ and moreover $H_{a,b}({\mathcal X}_a(0),{\mathcal X}_b(0)) = 1$.

b) Similarly, if $\varphi(b) \ne 0$ we had defined
\[
T_{a,b}(z,w) = \frac {w+1}{w} \left( 1 - \frac {z}{F_{a,b}(K_a(z),{\mathcal X}_b(w))}\right)
\]
where $(z,w) \in ({\mathbb C}\backslash\{0\})^2$ are close to $(0,0)$. Here 
\[
\begin{split}
F_{a,b}(t,s) &= \varphi((t1-a)^{-1}(1-sb)^{-1}) \\
&= \sum_{p \ge 0,q \ge 0} t^{-p-1}s^q\varphi(a^pb^q) \\
&= G_a(t) + s\theta(t,s)
\end{split}
\]
where $\theta(t,s) \in \sum_{p \ge 0,q \ge 0} t^{-p-1}s^{q-1}\varphi(a^pb^p)$ is holomorphic in a neighborhood of $(\infty,0) \in ({\mathbb C} \cup \{\infty\}) \times {\mathbb C}$. Since $G_a(K_a(z)) = z$, $K_a(z) = z^T{-1}{\tilde K}(z)$ where ${\tilde K}(z) = 1 + zR_a(z)$ is holomorphic near $0$ and ${\mathcal X}_b(s) = w{\tilde {\mathcal X}}(w)$ where ${\tilde {\mathcal X}}$ is holomorphic for $w$ in a neighborhood of $0$ we get that
\[
\begin{split}
T_{a,b}(z,w) &= \frac {w+1}{w} \frac {{\mathcal X}_b(w) \Theta (K_a(z),{\mathcal X}_b(w)}{z^{-1}F_{a,b}(K_a(z),{\mathcal X}_b(w))} \\
&= \frac {(w+1){\tilde {\mathcal X}}(w) \Theta (K_a(z),{\mathcal X}_b(w))}{z^{-1}F_{a,b}(K_a(z),{\mathcal X}_b(w))}\,.
\end{split}
\]
Since $K_a$ is analytic from a neighborhood of $0$ to a neighborhood of $\infty$, $\Theta(K_a(z),{\mathcal X}_b(w))$ is holomorphic in a neighborhood of $(0,0) \in {\mathbb C}^2$ and thus the numerator in the last formula is holomorphic near $0$. On the other hand the denominator is
\[
\sum_{p \ge 0,q \ge 0} z^p({\tilde K}(z))^{-p-1}({\mathcal X}_b(w))^q\varphi(a^pb^p)
\]
which is holomorphic in a neighborhood of $(0,0) \in {\mathbb C}^2$ and equal $1$ when $z = 0$, $w = 0$. Hence $T_{a,b}(z,w)$ extends to a holomorphic function of $(z,w)$ in a neighborhood of $(0,0)$.\qed

\bigskip
\noindent
{\bf Remark 4.1.} In view of Proposition~4.1 we will go beyond Definition~2.1 and Definition~3.1 from now on and view $S_{a,b}(z,w)$ and $T_{a,b}(z,w)$ as defined in a neighborhood of $(0,0)$.

\subsection{The case of factorizing $2$-band moments}
\label{sec4.3} In a sense the transforms $S_{a,b}(z,w)$ and $T_{a,b}(z,w)$, behave ``like covariances'', that is they are trivial in the case of $a$ and $b$ which behave like classically independent.

\bigskip
\noindent
{\bf Proposition 4.2.} {\em Let $a,b \in ({\mathcal A},\varphi)$ be such that $\varphi(a^pb^q) = \varphi(a^p)\varphi(a^q)$ for al $p \ge 0$, $q \ge 0$. Then if $\varphi(b) \ne 0$ we have $T_{a,b}(z,w) = 1$ and if additional also $\varphi(a) \ne 0$ then
\[
S_{a,b}(z,w) = 1.
\]
}

\bigskip
\noindent
{\bf {\em Proof.}} Under the assumptions in the proposition we have $F_{a,b}(t,s) = G_a(t)h_b(s)$ and $H_{a,b}(t,s) = h_a(t)h_b(s)$. Hence we have
\[
F_{a,b}(K_a(z),{\mathcal X}_b(w)) = G_a(K_a(z))H_a({\mathcal X}_b(w)) = z(w+1)
\]
and
\[
H_{a,b}({\mathcal X}_a(z),{\mathcal X}_b(w)) = h_a({\mathcal X}_a(z))h_b({\mathcal X}_b(w)) = (z+1)(w+1).
\]
This gives
\[
T_{a,b}(z,w) = \frac {w+1}{w} \left( 1 - \frac {z}{z(w+1)}\right) = 1
\]
and
\[
S_{a,b}(z,w) = \frac {z+1}{z} \frac {w+1}{w} \left( 1 - \frac {1+z+w}{(z+1)(w+1)} \right) = 1.
\]
\qed

\bigskip
\noindent
{\bf Remark 4.2.} Perhaps an adjective to point out this behavior of these transforms could be to say they are ``reduced'' transforms. For the partial two-variable $R$-transform there is a part of it which should be called the ``reduced'' transform, we will use ${\tilde R}$ instead of $R$ for it:
\[
{\tilde R}_{a,b}(z,w) = 1 - \frac {zw}{G_{a,b}(K_a(z),K_b(w))}\,.
\]
Under the factorization assumption in Proposition~4.2 we have $G_{a,b}(z,w) = G_a(z)G_b(w)$ and here
\[
G_{a,b}(K_a(z),K_b(w)) = G_a(K_a(z))G_b(K_b(w)) = zw
\]
so that ${\tilde R}_{a,b}(z,w) = 0$ in this case.

\subsection{Operations applied to probability measures yield signed measures}
\label{sec4.4} In free probability and free convolution operations yield free convolution operations on certain probability measures. The same happens with additive bi-free convolution and compact support probability measures on ${\mathbb R}^2$. It seems that for the operations $\boxtimes\boxtimes$ and $\boxtimes\boxplus$ and probability measures with compact support on $(0,\infty) \times (0,\infty)$ and, respectively, $(0,\infty) \times {\mathbb R}$ one gets only signed measures with finite total variation and integral $1$. This is easily seen from the following considerations. In the first case probability measures on $(0,\infty) \times (0,\infty)$ with compact support $\mu$ and $\nu$ will be the distribution of two two-faced pairs $(a_1,b_1)$ and $(a_2,b_2)$ in a $C^*$-probability space $({\mathcal A},\varphi)$ which are bi-free and where $a_k,b_k$ are $\ge 0$ and invertible and $[a_k,b_l] = 0$, $1 \le k$, $l \le 2$. Then the moments $\varphi((a_1a_2)^p(b_1b_2)^q)$ $p \ge 0$, $q \ge 0$ can also be written in the form $\varphi(a_1^{1/2}b_1^{1/2}(a_1^{1/2}a_2a_1^{1/2})^p(b_1^{1/2}b_2b_1^{1/2})^qb_1^{-1/2}a_1^{-1/2})$. Clearly $a_1^{1/2}a_2a_1^{1/2}$ and $b_1^{1/2}b_2b_1^{1/2}$ are commuting positive invertible operators, while the functional $\psi$ on ${\mathcal A}$ defined by $\psi(a) = \varphi(a_1^{1/2}b_1^{1/2}ab_1^{-1/2}a_1^{-1/2})$ is a bounded functional with $\psi(1) = 1$. This translates into $\mu\boxtimes\boxtimes\nu$ being a complex measure of bounded total variation and integral $1$. Similarly for $\boxtimes\boxplus$ we consider $a_k \ge 0$ invertible, $b_k = b_k^*$, $[a_k,b_l] = 0$ and $(a_1,b_1)$ and $(a_2,b_2)$ bi-free in $({\mathcal A},\varphi)$. Then we have that the moments $\varphi((a_1a_2)^p(b_1+b_2)^q)$ can also be written $\varphi(a_1^{1/2}(a_1^{1/2}a_2a_1^{1/2})^p(b_1+b_2)^qa_1^{-1/2})$ or $\psi((a_1^{1/2}a_2a_1^{1/2})^p(b_1+b_2)^q)$ where $\psi(a) = \varphi(a_1^{1/2}aa_1^{-1/2})$. Again the result $\mu\boxtimes\boxplus\nu$, where $\mu,\nu$ are probability measures on $(0,\infty) \times {\mathbb R}$, is a complex measure on $(0,\infty) \times {\mathbb R}$ with compact support, finite total variation and integral $1$. To see that the complex measures are actually signed measures, it suffices to remark that the moments $\varphi((a_1a_2)^p(b_1b_2)^q)$ in the first case and $\varphi((a_1a_2)^p(b_1+b_2)^q)$ in the second case are real numbers. This is indeed so for the general reason that $(a_k,b_k)$ can be realized, as long as we look only at moments, as the multiplication operators in the real Hilbert spaces $L_{\text{red}}^2({\mathbb R}^2,d\mu)$, $L_{\text{real}}^2({\mathbb R}^2,d\nu)$ with state vector the constant function $1$. The free product of these real vector spaces with state-vector can be then carried out over ${\mathbb R}$ and the left and right operators we get on the real free product have then real moments.

\subsection{The $S$- and $T$-bifree partial transforms in the general algebraic setting and formal power series}
\label{sec4.5} The two-variables partial $S$- and $T$-transforms can also be defined in the purely algebraic context as formal power series with complex coefficients and the analogues of Theorem~2.1 and Theorem~3.1 remain valid. This is due to the fact that these relations boil down to polynomial relations between moments under bi-freeness assumptions an that the non-commutative variables in the Banach-algebra context are sufficiently general to insure that the polynomial relations hold in general. (Concerning the free convolutions handled in very general algebraic settings see \cite{4}.)

\bigskip
\noindent
{\bf Note added May 14, 2015:} Soon after this paper appeared as preprint arXiv:1504.03765, combinatorial proofs for the $S$- and $T$-transforms were found by Paul Skoufranis in his preprint ``A combinatorial approach to Voiculescu's bi-free partial transforms'' arXiv:1504.06005.

\end{document}